\documentclass[12pt]{article}
\usepackage[margin=1in]{geometry}
\usepackage{amsmath,amssymb,amsthm,amsfonts,color,hyperref,graphicx,subfigure,authblk}
\usepackage{indentfirst}
\hypersetup{colorlinks,
linkcolor=blue} 

\theoremstyle{definition}
\newtheorem{thm}{Theorem}[section]
\newtheorem{lem}[thm]{Lemma}

\newtheorem{prop}[thm]{Proposition}
\newtheorem{cor}[thm]{Corollary}

\numberwithin{equation}{section}

\newcommand{\enabstractname}{Abstract}
\newenvironment{enabstract}{
    \par\small
    \noindent\mbox{}\hfill{\bfseries \enabstractname}\hfill\mbox{}\par
    \vskip 2.5ex}{\par\vskip 2.5ex} 

\title{A N-body problem with weak force potential through Hamilton-Jacobi equation approach}
\author[a]{Putian Yang\thanks{Corresponding author: yangputian@stu.scu.edu.cn}}
\author[a]{Shiqing Zhang}

\affil[a]{Department of mathmatics, Sichuan University}

\date{\today}

\begin{document}

\maketitle

\begin{enabstract}
    This paper we consider for the $N$-body problem with potential $1/r^{\alpha}$
    ($0<\alpha<1$)
     the existence of hyperbolic motions for any prescribed 
    limit shape and any given initial configuration of the bodies. 
    Here $E$ is the Euclidean space where the bodies moving and $\|\cdot\|_{E}$ is 
    the norm induced by the inner product.
     The energy level $h>0$ of 
    the motion can also be chosen arbitrarily. We use the
    global viscosity solutions for the Hamilton-Jacobi equation $H\left(x, d_{x} u\right)=h$ and geodesics.

\end{enabstract}

\tableofcontents

\newpage

\section{Introduction}
  Before we discribe the background and main issue of this paper, we firstly present some of notations about our working space $E=\mathbb{R}^{n}, (n\geq 2)$.
    $\langle\cdot,\cdot \rangle $ is the standard inner product of $ E$,
  $|\cdot|$ is its induced norm. We also write the inner product of configuration space $E^{N}$ as 
  $$
  \langle x,y\rangle=\left(\frac{1}{2}\sum_{i=1}^{N}\sqrt{m_{i}m^{'}_{i}}\langle x_{i},y_{i}\rangle\right)^{1/2}
  $$ for $x=(x_{1},\ldots,x_{N}), y=(y_{1},\ldots,y_{N})\in E^{N}$ where $m_{i}$ and $m^{'}_{i}$ are mass of  $i$-th body of configurations $x$ and $y$ respectively,
   and  $x_{i},y_{i}\in E,(1\leq i\leq N)$. 
   The induced norm of configuration space $E^{N}$ is denoted as $\|\cdot\|$ and defined as $\|x\|=\left(\frac{1}{2}\sum_{i=1}^{N}m_{i}|x_{i}|^{2}\right)^{1/2}.$
  Write $S_{\lambda}:=\{(x,y)\in E\times E\mid|x-y|<\lambda\}$.

  We study N-body problem in $ E$ with potential
   $$U(x)=\sum_{1\leq i<j\leq N}\frac{m_{i}m_{j}}{|x_{i}-x_{j}|^{\alpha}}$$
   where $x=(x_{1},\ldots,x_{N})\in E^{N}, (x_{i}\in  E)$ is the configuration, $m_{i}$ is the mass of $i$-th body of $x$.

Our main task is to find the hyperbolic motions of our stated N-body problem by measure of construction, and the motion satiesfies
$$
\ddot{x_{i}}=\nabla_{x_{i}}U(x)=\alpha\sum_{j=1,j\neq i}^{N}\frac{x_{j}-x_{i}}{|x_{j}-x_{i}|^{\alpha+2}}.
$$
i.e., 
\begin{equation}\label{1.1}
  \ddot{x}=\nabla U(x)
\end{equation}

Now we give some notations frequently used in this paper.
\begin{enumerate}
  \item Denote $\Omega$ be the subset of configurations with non-collision, i.e., 
$$
\Omega=\left\{x\in E^{N}\mid x_{i}\neq x_{j},  \forall 1\leq i<j\leq N\right\}
.$$
write $\sum:=E^{N}\setminus\Omega$.

\item Denote $r(x)=\min \{|x_{i}-x_{j}|\}_{1\leq i<j\leq N}$, this is a essential data for the judgement of collision of a configuration is the minimal distance among its bodies, meaning that
 for $x\in E^{N}$, $r(x)>0$ if and only if $x\in\Omega$. Similarly we denote $R(x)=\max \{|x_{i}-x_{j}|\}_{1\leq i<j\leq N}$.
 \item Denote several line segement as follows:

  \begin{equation}
    \overline{[xy]}=\left\{tx+(1-t)y\mid x,y\in E^{N}; 0\leq t\leq 1 \right\}
  \end{equation} 

  \begin{equation}
    \overline{(xy)}=\left\{tx+(1-t)y\mid x,y\in E^{N}; 0< t< 1 \right\}
  \end{equation} 

  \begin{equation}
    \overline{(xy\infty)}=\left\{x+ty\mid x,y\in E^{N}; t>0\right\}
  \end{equation} 

  \begin{equation}
    \overline{[xy\infty]}=\left\{x+ty\mid x,y\in E^{N}; 0\leq t\leq+\infty \right\}.
  \end{equation}

\item Denote $d(x,W)=\inf\left\{\|x-y\|\mid y\in W\right\}$  as
 the distance of a configuration $x\in E^{N}$ to some subset $W\subset E^{N}$.

\item Denote
\begin{equation}
\mathbf{C}(x,y,T):=\left\{\gamma:[\alpha,\beta]\rightarrow E^{N}\text{ for some }[\alpha,\beta]\subset \mathbb{R}\mid\gamma(\alpha)=x,\gamma(\beta)=y,\beta-\alpha=T\right\}
\end{equation}
and
\begin{equation}
\mathbf{C}(x,y)=\cup_{T>0} \mathbf{C}(x,y,T),
\end{equation}
without loss of generality, we can always assume $\alpha=0, \beta=T>0$ for $T$ under determinded.

\item Denote \begin{equation}
  l(\gamma):=\int_{0}^{T}|\dot{\gamma}(t)|dt 
\end{equation}
as the Euclidean length of $\gamma\in \mathbf{C}(x,y)$.

\item Denote $\angle(x,y)$ as the angle between $x,y\in E$, i.e., 
\begin{equation}
  \cos\angle(x,y):=\frac{\langle x,y\rangle}{|x||y|}.
\end{equation}
Similarly, if $x,y\in E^{N}$,
\begin{equation}
  \cos\angle(x,y):=\frac{\langle x,y\rangle}{\|x\|\|y\|}.
\end{equation}
\item Denote $\mathbb{S}^{nN-1}=\{x\in E^{N}\mid \|x\|=1 \}$ where $n=\dim E$, we find that for $x,y\in\mathbb{S}^{nN-1}$, 
\begin{equation}
  \frac{1}{2}\|x-y\|=\sin\frac{1}{2}\angle(x,y)
\end{equation}
\item

\end{enumerate}

We give several assumption that is for simplicity thourghout the artical.
\begin{enumerate}
  \item We assume $\min\{m_{i}\}_{1\leq i\leq N}=1$.
\end{enumerate}

Now we can present the main thereom as our main issue 
\begin{thm}
  There is a hyperbolic motion $\gamma=(\gamma_{1},\ldots,\gamma_{N}):[0,+\infty)\rightarrow E^{N}$ which $\gamma(0)=x_{0}$ and $\gamma(t)=\sqrt{2E}at+o(t)$ as $t\rightarrow+\infty$ where\\
  1) $x_{0}\in E^{N}$ is arbitrarily given initial configuration.                       \\
  2) $E>0$ is arbitrarily given energy constant.                      \\
  3) $a\in E^{N}$ is arbitrarily given non-collision configuration and $\|a\|=1$

  In particular, $\gamma\subset\Omega$ and is a solution of \ref{1.1}.
\end{thm}

We check $TE^{N}=E^{N}\times(E^{N})^{\star}\simeq E^{N}\times E^{N}$ and the Hamiltonian on $TE^{N}=E^{N}\times(E^{N})^{\star}$:  $H(x,p)=\|p\|^{2}-U(x), p\in E^{N}.$

and the Lagrangian $L$ on $E^{N}\times E^{N}$ as the 
dual of the Hamiltonian:
\begin{equation}
  L(x,v)=\sup\left\{p(v)-H(x,p)\mid p\in E^{N}\right\}=\|v\|^{2}+U(x)
\end{equation}
and 
\begin{equation}
  H(x,p)=\sup\left\{p(v)-L(x,v)\mid v\in E^{N}\right\}
\end{equation}
here $L(x,v)=+\infty$ when $x$ is a collision.

We define the Lagrangian action of $\gamma\in\mathbf{C}(x,y)$ as 
$$
\mathbf{A}_{E}(\gamma)=\int_{\alpha}^{\beta}L(\gamma,\dot{\gamma})+Edt=\int_{\alpha}^{\beta}\|\dot{\gamma}\|^{2}+U(\gamma)+Edt
,$$ where $E>0$ is the energy constant. 
We notice that this is well defined since $\Omega$ is connected, for all $x,y\in\Omega$ there is a smooth $\gamma\in\mathbf{C}(x,y,T)$ and $U$ is bounded 
in the compact set $\gamma|_{[0,T]}$, therefore $\mathbf{A}_{E}(\gamma)<\infty$.
Thus we can define the minimal action function $\phi_{E}(\cdot,\cdot)$ as:
$$
\phi_{E}(x,y)=\inf\{\mathbf{A}_{E}(\gamma)\mid\gamma\in\mathbf{C}(x,y)\}.
$$
The specific $\gamma\in\mathbf{C}(x,y)$ which satiesfies $\phi_{E}(x,y)=\mathbf{A}_{E}(\gamma)$ is called a free time minimizer of the action $A_{E}$.

Maderna presented in his another paper\cite{da2014free} the existence of free time minimizers in the Newtonian N-body problem, on the otherhand, 
we here present a work by Mather\cite{mather1991action} a version of Tonelli's theorem which we state below,
for a simplified version see \cite{contreras1999global}.
\begin{thm}
  For all $x, y\in E^{N}$ there exists a Tonelli minimizer on $\mathbf{C}(x,y)$, in the sense that
  $A_{E}(\gamma)=\min_{\eta \in \mathbf{C}(x, y)} A_{E}(\eta)$
\end{thm}

The singularity problem of motions of N-body is well known and widely studied, for example,
we already know that in 2002 Marchal presented in\cite{marchal2002method} that the minimization process confirms the absence of collision in the classical Newtonian N-body problem.
\begin{thm}(2002, Marchal \cite{marchal2002method}). 
  If $\gamma \in \mathbf{C}(x, y)$ is defined on some interval $[a, b]$, and satisfies $\mathbf{A}_{L}(\gamma)=\phi(x, y, b-a)$, then $\gamma(t) \in \Omega$ for all $t \in(a, b)$.
\end{thm}

Here in our paper we need to know the absence of collision of potential with homogeneity $-\alpha<0$, they are studied, see[]
for the convienence we restate here.
\begin{thm}
  When $U(x)=\sum_{1\leq i<j\leq N}\frac{m_{i}m_{j}}{|x_{i}-x_{j}|^{\alpha}}$
  All minimizers of $\phi_{E}$ are experence no collison, i.e., if $A_{E}(\gamma)=\phi_{E}(x.y)$ where $\gamma\in\mathbb{A}(x,y,T)$ for
  some $T>0$, then $\gamma(t)\in\Omega, t\in(0,T)$
\end{thm}

We also have to list the very important Hamilton's principle of least action as the following thereom which we will use and its proof is put later.
\begin{thm}
  For any $E>0$, suppose $\gamma\in\mathbf{C}(x,y,T)$ ($T$ can be $+\infty$)  is a free time minimizer of $A_{E}$ and $\gamma|_{(0,T)}\subset\Omega$,
  then $\gamma$ is a solution of \ref{1.1} in $(0,T)$.
 \end{thm}

There is another result of $\phi_{E}(\cdot,\cdot)$.
\begin{prop}
  $\phi_{E}(\cdot,\cdot)$ is a distence function in $\Omega$
\end{prop}
\begin{proof}
  We first notice $\phi_{E}$ meets trianglular innequality.

  For any $x,y,z\in\Omega$, any $\gamma_{1},\gamma_{2}\in\mathbf{C}(x,y)$ defined in $[0,T_{1}]$ and $[0,T_{2}]$ respectively, we set 
\begin{equation}
  \gamma(t):=
  \begin{cases}
    \gamma_{1}(t) &0\leq t\leq T_{1}\\
    \gamma_{2}(t-T_{1}) &T_{1}\leq t\leq T_{1}+T_{2}
  \end{cases},
\end{equation}
thus \begin{equation}
  \phi_{E}(x,z)=\inf\left\{\mathbf{A}_{E}(\eta)\mid\eta\in\mathbf{C}(x,z)\right\}\leq\mathbf{A}_{E}(\gamma)\leq\mathbf{A}_{E}(\gamma_{1})+\mathbf{A}_{E}(\gamma_{2}).
\end{equation}
Since $\gamma_{1},\gamma_{2}$ are arbitrary, we have $\phi_{E}(x,z)\leq\phi_{E}(x,y)+\phi_{E}(y,z)$.

Second, we verify that $\phi_{E}(x,y)=0$ makes $x=y$.

For any $x=(x_{1},\ldots,x_{N}),y=(y_{1},\ldots,y_{N})\in E^{N}$ and $\gamma=(\gamma_{1},\ldots,\gamma_{N})\in\mathbf{C}(x,y)$ defined on $[0,T]$.
There is a $T'\in(0,T]$ s.t. $$\max\left\{|\gamma_{i}(T')-x_{i}|\mid 1\leq i\leq N\right\}=\max\left\{|x_{i}-y_{i}|\mid 1\leq i\leq N\right\},$$
and $$\max\left\{|\gamma_{i}(t)-x_{i}|\mid 1\leq i\leq N\right\}\leq\max\left\{|x_{i}-y_{i}|\mid 1\leq i\leq N\right\}$$ for $t\in[0,T')$.
Then there exists $i_{0}$ s.t. 
\begin{equation} 
  \max\left\{|x_{i}-y_{i}|\mid 1\leq i\leq N\right\}=\left|\gamma_{i_{0}}(T')-x_{i_{0}}\right|\leq\int_{0}^{T'}|\dot{\gamma_{i_{0}}}|dt\leq \sqrt{T'}\left(\int_{0}^{T'}|\dot{\gamma_{i_{0}}}|^{2}dt\right)^{1/2}.
\end{equation}
Hence  
\begin{align}
  \mathbf{A}_{E}(\gamma)&\geq \mathbf{A}_{E}\left(\gamma\left|_{[0,T']}\right.\right)\geq\int_{0}^{T'}\|\dot{\gamma}(t)\|^{2}+U\left(\gamma(t)\right)+Edt\\
                        &\geq \int_{0}^{T'}\|\dot{\gamma}(t)\|^{2}dt=\frac{1}{2}\int_{0}^{T'}\sum m_{i}|\dot{\gamma}_{i}(t)|^{2}dt\geq \frac{m_{i_{0}}}{2}\int_{0}^{T'}|\dot{\gamma}_{i_{0}}(t)|^{2}dt\\
                        &\geq \frac{1}{2T'}\max\left\{|x_{i}-y_{i}|^{2}\mid 1\leq i\leq N\right\}.
                      \end{align}
So eventually $\phi_{E}(x,y)=0$ makes $\max\left\{|x_{i}-y_{i}|\right\}_{i=1}^{N}=0$, which means $x=y$.

It is not difficult to see $\phi_{E}(x,x)=0$.
\end{proof}

\section{The existence and properties of free-time minimizers}

\section{Some preparations for geometric objects}
Based on the elementary computing, we can verify the following essential geometric facts in Euclidean space $E^{N}$.

\begin{lem}\label{lem1.1}
  For any $x\in E^{N}$, we must have $|x_{i}|\leq\sqrt{2}\|x\|$ and $|x_{i}-x_{j}|<3\|x\|$ for $1\leq i<j\leq N$, hence $r(x)<3\|x\|$.
  \end{lem}
\begin{proof}
  By the definition of $\|x\|$, we have $\sum|x_{i}|^{2}\leq\sum m_{i}|x_{i}|^{2}=2\|x\|^{2}$, thus $|x_{i}|\leq\sqrt{2}\|x\|$ for any $i$ and $|x_{i}-x_{j}|\leq 2\sqrt{2}\|x\|<3\|x\|$.

\end{proof}

\begin{thm}
  For $a\in\mathbb{S}^{nN-1}, r(a)>0, x,x'\in E^{N}$, the following statements are valid.
  \begin{enumerate}
    \item \begin{equation}
      \left\lVert \frac{x+ta}{\|x+ta\|}-a \right\rVert \leq\frac{1}{30}r(a)
    \end{equation}
    and 
    \begin{equation}
      r(x+ta)\geq\frac{67}{70}r(a)t>67
    \end{equation}
    for any $t>\frac{1+\|x\|}{r(a)}$

    \item
     \begin{equation}
      r(x')>r(a)-2\|x'-a\|\geq(1-3\lambda)r(a)
    \end{equation}
    and
    \begin{equation}
      \cos\angle(x'_{i}-x'_{j}, a_{i}-a_{j})\geq 1-6\lambda, 1\leq i<j\leq N
    \end{equation}
    
    for $\|x'-a\|\leq\lambda r(a), \lambda\in(0,1/2)$
    
\item If $\|x'\|=1$ then
\begin{equation}
  \cos\angle(a,x')=\langle a,x'\rangle\geq 1-\frac{9}{2}\lambda^{2}
\end{equation}
  \end{enumerate} 
  
\end{thm} 
\begin{cor}
\end{cor}

\begin{proof}
  Since $t>70\frac{1+\|x\|}{r(a)}$, we have $t>70\frac{\|x\|}{r(a)}$,  $\|x\|<\frac{r(a)}{70}t$ and $r(a)t>70$.

   \begin{enumerate}
    \item Since $t>70\frac{\|x\|}{r(a)}$ and by lemma \ref{1.1} we know that $r(a)<3\|a\|=3$, we have $t-\|x\|>\frac{70-3}{r(a)}\|x\|>\frac{60}{r(a)}\|x\|$ hence $\frac{\|x\|}{t-\|x\|}<\frac{r(a)}{60}$.
  \begin{align}
    \left\lVert \frac{x+ta}{\|x+ta\|}-a\right\rVert&=\frac{1}{\|x-ta\|}\left\lVert x+ta-\|x+ta\|a\right\lVert\\
    &\leq\frac{1}{\|x+ta\|}\left(\|x\|+\left\lVert \|ta\|-\|x+ta\| \right\rVert \right)\leq\frac{2\|x\|}{\|x+ta\|}\\
    &\leq\frac{2\|x\|}{t-\|x\|}\leq\frac{r(a)}{30}.
  \end{align}

\begin{align}
  |x_{i}+ta_{i}-(x_{j}+ta_{j})|&=|t(a_{i}-a_{j})+x_{i}-x_{j}|\\
  &\geq t|a_{i}-a_{j}|-|x_{i}-x_{j}|\\
  &\geq r(a)t-3\|x\|>r(a)t-2\frac{r(a)}{70}t=\frac{67}{70}r(a)t>67.
\end{align}
\item First we have 
\begin{equation}
  |a_{i}-a_{j}|\leq|a_{i}-x'_{i}|+|x'_{i}-x'_{j}|+|x'_{j}-a_{j}|.
\end{equation}
thus 
\begin{equation}
  |x'_{i}-x'_{j}|\geq|a_{i}-a_{j}|-|x'_{i}-a_{i}|-|x'_{j}-a_{j}|>|a_{i}-a_{j}|-3\|x'-a\|,
\end{equation}
hence 
\begin{equation}
  r(x')>r(a)-3\|x'-a\|>(1-3\lambda)r(a).
\end{equation}
On the otherhand 
\begin{align}
  \langle a_{i}-a_{j}, x'_{i}-x'_{j} \rangle &= \langle a_{i}-a_{j}, a_{i}-a_{j} \rangle + \langle a_{i}-a_{j}, x'_{i}-a_{i}-x'_{j}+a_{j} \rangle\\
  &=|a_{i}-a_{j}|^{2}+\langle a_{i}-a_{j}, (x'_{i}-a_{i})-(x'_{j}-a_{j})\rangle\\
  &\geq |a_{i}-a_{j}|^{2}-|a_{i}-a_{j}||(x'_{i}-a_{i})-(x'_{j}-a_{j})|\\
  &\geq |a_{i}-a_{j}|^{2}-|a_{i}-a_{j}|3\|x'-a\|\\
  &=|a_{i}-a_{j}|\left(|a_{i}-a_{j}|-3\|x'-a\|\right).
\end{align}

\begin{equation}
  |x'_{i}-x'_{j}|\leq|x'_{i}-a_{i}|+|a_{i}-a_{j}|+|a_{j}-x'_{j}|\leq|a_{i}-a_{j}|+3\|x'-a\|.
\end{equation}
So we eventually have 
\begin{align}
  \frac{\langle a_{i}-a_{j}, x'_{i}-x'_{j}\rangle}{|a_{i}-a_{j}||x'_{i}-x'_{j}|}&\geq \frac{|a_{i}-a_{j}|-3\|x'-a\|}{|a_{i}-a_{j}|+3\|x'-a\|}\\
  &\geq 1-\frac{6\|x'-a\|}{|a_{i}-a_{j}|+3\|x'-a\|}\\
  &\geq 1-\frac{6\|x'-a\|}{r(a)}\geq 1-6\lambda.
\end{align}
\item
For $\|x'\|=1$ 
\begin{equation}
  2-2\langle x',a\rangle=\|x'-a\|^{2}\leq\lambda^{2}r(a)^{2}.
\end{equation}
hence 
\begin{equation}
  \langle x',a\rangle=1-\frac{1}{2}\|x'-a\|^{2}\geq 1-\frac{\lambda^{2}}{2}r(a)^{2}=1-\frac{9}{2}\lambda^{2}
\end{equation}
  \end{enumerate}
\end{proof}

\bibliographystyle{plain}
\bibliography{ref}
\addcontentsline{toc}{section}{Reference}

\end{document}